\documentclass{article}
\usepackage{latexsym}
\usepackage{amssymb}
\usepackage[dvips]{epsfig}

\def\Z{{\mathbf Z}}

\def\R{{\mathbf R}}
\def\C{{\mathbf C}}
\def\T{{\mathbf T}}

\def\cI{{\cal I}}

\def\i{{\imath}}

\def\cInf{\mathop{C^{\infty}}\nolimits}

\def\span{\mathop{\rm span}\nolimits}

\def\Id{\mathrm{id}}

\def\Im{\mathrm {Im}}
\def\Re{\mathrm {Re}}
\def\Spec{\mathrm{Spec}}

\def\Homsg{\mathrm{Hom_{s.g.}}}
\def\Hom{\mathrm{Hom}}
\def\Int{\mathrm{int}}

\def\Box{\hspace{10pt}\mathfrak {qed}}
\def\:{\stackrel{\small\mathrm{def}}{\small{=}}}

\newcommand{\Pbr}[2]{\{#1,#2\}}

\newcommand{\Eval}[2]{\langle #1,#2\rangle}

\newcommand{\D}[1]{\frac{\partial}{\partial #1}}
\newcommand{\p}[2]{\frac{\partial #1}{\partial #2}}

\newcommand{\DD}[2]{\frac{\partial}{\partial#1}\wedge\frac{\partial}{\partial#2}}

\newcommand{\cc}[1]{\overline #1}
\newcommand{\n}[2]{{#1}^{\ast #2}}
\newcommand{\iw}[2]{#1_{#2_1, \ldots, #2_k}}
\newcommand{\m}[1]{\mathfrak #1}
\newcommand{\binom}[2]{{#1 \choose #2}}
\newtheorem{theorem}{Theorem}

\newtheorem{lemma}{Lemma}
\newtheorem{corollary}{Corollary}

\newtheorem{example}{Example}

\title{Newton polyhedra and Poisson structures from certain linear Hamiltonian circle actions}

\author{\'{A}g\'{u}st Sverrir Egilsson\\
Department of Mathematics,\\
University of California, Berkeley.\\
egilsson@math.berkeley.edu}

\begin{document}
\maketitle

\begin{abstract}
In this paper we first describe the geometry of the Newton polyhedra
of polynomials invariant under certain linear Hamiltonian circle
actions. From the geometry of the polyhedra, various Poisson
structures on the orbit spaces of the actions are derived and
Poisson embeddings into model spaces, for the orbit spaces, are
constructed. The Poisson structures, on respective source and model
space, are compatible even for the minimum possible (embedding)
dimension of the model spaces. This is, in particular, important
since it is still an open question if, in general, there exist
finite dimensional model spaces with Poisson structures compatible
with the actions and the usual nondegenerate Poisson structure on
the source spaces.
\end{abstract}
\tableofcontents

\section{Prefatory notes} To understand the motivation for
the subject at hand, it may be enlightening to visit the 3 examples
in section 8.1 of \cite{ASE1}, as well as the main result from that
earlier paper: Namely, the rarity of the existence of Poisson
structures, on the model spaces, compatible with the usual Poisson
structure on the orbit spaces. In the paper at hand, we relax the
strict requirements, used in \cite{ASE1}, of starting with
symplectic source spaces and instead equip the spaces with Poisson
structures derived from the combinatorics observed by applying the
Hamiltonian actions.

\section{Introduction}

In this introductory section of the paper, simple results obtained
for the circle action
$$
\T\times \C^{k} \rightarrow \C^{k}\mbox{ where }z (z_{1}, \ldots,
z_{k}) = (z^{n_{1}}z_{1}, \ldots, z^{n_{k}}z_{k}),
$$
where $\T$ is the unit circle in $\C$ and $n_1,\ldots,n_k$ are
nonzero integers, called weights, are recalled. A more complete
discussion is available in \cite{ASE1} \S 1-\S 3. A polynomial in
$\C[z_1,\ldots,z_k,\cc{z}_1,\ldots,\cc{z}_k]$ is invariant under the
action if and only if the exponents of each of its nonzero terms
$z_1^{a_1} \cc{z}_1^{b_1} \cdots z_k^{a_k} \cc{z}_k^{b_k}$ satisfies
the equation $\n{n}{}(a-b) ~[= (a_1-b_1) n_1 + \cdots + (a_k-b_k)
n_k] = 0$. The monoid of all such lattice points $(a,b)$ determining
the invariant polynomials for a fixed set of weights
$n_1,\ldots,n_k$ is denoted by $(\iw{S}{n},0,+)$, the ring of
invariant polynomials is $\C[\iw{S}{n}]$ and the group generated by
$\iw{S}{n}$ in $\Z^k \times \Z^k$ is denoted by $\iw{M}{n}$. It,
i.e., $\iw{M}{n}$, consists of all solutions $(a,b)\in \Z^k \times
\Z^k$ to $\n{n}{}(a-b)=0$. The usual basis for $\Z^k \times \Z^k$ is
denoted by the elements $e_1,\ldots,e_k,\cc{e}_1,\ldots,\cc{e}_k$.

Restricting to multiple [$k>1$] positive relative prime
[$\gcd(n_1,\ldots,n_k)=1$] weights $n_1,\ldots,n_k$, let $\iota$ be
the map from the integer hyperplane $n^\perp = \{r\in \Z^k : r_1
n_1+ \cdots + r_k n_k = 0\}$ to the integer hyperplane $\cI^\perp =
\{t\in \Z^k : t_1 + \cdots + t_k = 0\}$ given by
$$
\iota: n^\perp \rightarrow \cI^\perp ; r \mapsto \frac{(n_1 r_1,
\ldots, n_k r_k)}{d_1 \cdots d_k}
$$
where each $d_i$ is defined by $d_i =
\gcd(n_1,\ldots,n_{i-1},n_{i+1},\ldots,n_k)$. The following
equivalences [i - iii] are observed in \cite{ASE1}: i) $\iota$ is an
isomorphism, ii) $n_i =  d_1 \cdots d_{i-1} d_{i+1} \cdots d_k$ for
$i = 1,\ldots,k$ and iii) the semigroup $\iw{S}{n}$ can be generated
by exactly $k^2$ elements. If these conditions hold then the
semigroup $\iw{S}{n}$ may be generated by the $k^2$ elements $e_1 +
\cc{e}_1, \ldots, e_k + \cc{e}_k$ and $d_i e_i + d_j \cc{e}_j$ for
$i\neq j$. To distinguish circle actions, or integers $n_1,\ldots,
n_k$, that satisfy these conditions we designate them as actions, or
integers, generating a minimal Hilbert basis [for the polynomials
invariant under the action].

The group $\iw{M}{n}$ is a lattice of rank $2k-1$ and has a basis of
$2k-1$ elements as follows.
\begin{lemma}\label{le-1}
Assume that $\iota$ is an isomorphism. Define elements
$l_1,\ldots,l_k$ by $l_i = e_i + \cc{e}_i$ and
$\eta_1,\ldots,\eta_{k-1}$ by $\eta_i = d_{i+1} e_{i+1} + d_i
\cc{e}_i.$  Then $l_1, \ldots, l_k, \eta_1,\ldots,\eta_{k-1}$ is a
basis for $\iw{M}{n}$.
\end{lemma}
{\it Proof: } Notice first that $\cI^\perp$ has a $\Z$-basis $\{e_2
- e_1,\ldots, e_k-e_{k-1}\}$ so $n^\perp$ has a basis containing the
$k-1$ elements $\iota^{-1}(e_{i+1}-e_i) = d_{i+1} e_{i+1} - d_i e_i
= \eta_i - d_i l_i$ for $i=1,\ldots,k-1$. Take $x=(a,b)\in
\iw{M}{n}$. Let $L = b_1 l_1 + \cdots + b_k l_k$, then $x - L \in
n^\perp = \span_\Z\{\eta_1 - d_1 l_1,\ldots,\eta_{k-1} - d_{k-1}
l_{k-1}\}$. Conclude that $x$ is in the $\Z$ span of the independent
elements $l_1, \ldots,l_k, \eta_1, \ldots, \eta_{k-1}$.
$\Box$\medskip

\section{Structure of the Newton polyhedra}

The finitely generated semigroup $\iw{S}{n}$ spans, over $\R_{\geq
0}$, a strongly convex rational polyhedral cone $\R_{\geq 0}\cdot
\iw{S}{n}$ in the $2k-1$ dimensional vector space $\R\cdot
\iw{M}{n}$. It is a convex polyhedral cone since it is the
non-negative span of finitely many vectors in a vector space,
strongly convex since it does not contain a line through the origin
and rational because it is generated by elements in the lattice
$\iw{M}{n}$. Assuming that the circle action generates a minimal
Hilbert basis, the extreme rays [one dimensional faces or edges] of
the cone are the $k^2$ half-lines $\R_{\geq 0} (d_i e_i + d_j
\cc{e}_j)$.

Let $\m{h}$ and $\m{v}$ be subsets of $\{1,\ldots,k\}$ and define
$\m{F}_{\m{h}\times\m{v}}$ to be the sub-cone of $\R_{\geq 0}\cdot
\iw{S}{n}$ given by the span
\[\m{F}_{\m{h}\times\m{v}} = \R_{\geq 0}\cdot\{ d_i e_i + d_j \cc{e}_j :
(i,j)\in\m{h}\times\m{v}\}.\] We now prove the following.

\begin{theorem}\label{th-1}
Assume that the integers $n_1,\ldots,n_k$ generate a minimal Hilbert
basis. Then the map
\[\m{h}\times\m{v} \rightarrow \m{F}_{\m{h} \times \m{v}}\] defines
a one-to-one relationship between subsets of
$\{1,\ldots,k\}\times\{1,\ldots,k\}$ of the form $\m{h}\times\m{v}$
and the faces of $\R_{\geq 0}\cdot \iw{S}{n}$.
\end{theorem}
In order to show this we use the following:
\begin{lemma}\label{le-2}
A strongly convex polyhedral cone $\tau$ contained as a subset in a
strongly convex polyhedral cone $\gamma$ is a face of $\gamma$ if
and only if for any $x$ and $y$ in $\gamma$ with the sum $x+y$ in
$\tau$ we have that $x\in \tau$ and $y\in\tau$.
\end{lemma}
{\it Proof of lemma: } First we assume that $\tau$ is a face of
$\gamma$, say $\tau = \gamma\cap v^{\perp}$ for some linear
functional nonnegative on $\gamma$. If $x+y\in\tau$ then
$\Eval{v}{x+y}=0$ but also $\Eval{v}{x}\geq 0$ and $\Eval{v}{y}\geq
0$ so we must have $\Eval{v}{x}=\Eval{v}{y}=0$.

\noindent Now, for the opposite direction, assume that for any
$x,y\in\gamma$ with $x+y\in \tau$ we have that $x,y\in\tau$. Let
$\tau'$ be a minimal face of $\gamma$ containing $\tau$. Let
$n=\dim(\tau')$. Assume $z\in\tau\cap\Int(\tau')$ [$\Int(\tau')$
denotes the relative topological interior of $\tau'$]. For any $x$
in a small $n$-ball in $\tau'$ centered at $z$ we can find $y$ in
the ball such that $z = x/2 + y/2$. Since $x/2$ and $y/2$ are in
$\gamma$ we have that $x/2$, and therefore also $x$, is in $\tau$.
We conclude that $\tau\cap\Int(\tau')$ is an open subset of
$\Int(\tau')$, it is also closed in $\Int(\tau')$ since $\tau$ is
closed. The set $\tau\cap\Int(\tau')$ is nonempty by the minimal
condition on $\tau'$. Since $\Int(\tau')$ is connected we obtain
that $\Int(\tau')\subset\tau\subset\tau'$. Hence $\tau' = \tau$, and
in particular $\tau$ is a face of $\gamma$.$\Box$\medskip

\noindent{\it Proof of theorem: } Using Lemma \ref{le-2} it follows
that each of the sets $\m{F}_{\m{h} \times \m{v}}$ is a face of
$\R_{\geq 0}\cdot \iw{S}{n}$. The map ${h}\times\m{v} \rightarrow
\m{F}_{\m{h} \times \m{v}}$ is clearly injective. Assume that $\tau$
is a face of $\R_{\geq 0}\cdot \iw{S}{n}$. We can assume that $\tau$
is spanned, over $\R_{\geq 0}$, by some extreme rays. Therefore,
$\tau$ has the format $\tau = \R_{\geq 0}\cdot\{d_i e_i + d_j
\cc{e}_j: (i,j) \in \m{i}\}$ for a set of indices $\m{i} \subset
\{1,\ldots,k\}\times\{1,\ldots,k\}$. For elements $(i',j')$ and
$(i'',j'')$ in $\m{i}$ we have that $(d_{i'} e_{i'} + d_{j'}
\cc{e}_{j'}) + (d_{i''} e_{i''} + d_{j''} \cc{e}_{j''}) = (d_{i'}
e_{i'} + d_{j''} \cc{e}_{j''}) + (d_{i''} e_{i''} + d_{j'}
\cc{e}_{j'})$ so by Lemma \ref{le-2} we also must have $(i',j''),
(i'',j') \in \m{i}$. From this it follows that $\m{i}$ has the cross
product format $\m{h}\times\m{v}$ for some
$\m{h},\m{v}\subset\{1,\ldots,k\}$ and therefore that $\tau =
\m{F}_{\m{h}\times\m{v}}$.  $\Box$\medskip

\noindent For a subset $\m{i}$ of the finite $k$-lattice,
$\m{i}\subset \{1, \ldots, k\} \times \{1, \ldots,k \}$, we define
\[\m{F}_{\m{i}} = \m{F}_{\m{h}\times\m{v}}\] where $\m{h}$ and
$\m{v}$ are the smallest sets such that
$\m{i}\subset\m{h}\times\m{v}$. We refer to elements of $\m{h}$ as
horizontal lines and elements of $\m{v}$ as vertical lines.
\begin{lemma}\label{le-3}
Assume that the integers $n_1,\ldots,n_k$ generate a minimal Hilbert
basis. Let $i$ be a non-empty subset of the finite $k$-lattice and
let $\m{h}$ and $\m{v}$ be minimal with respect to the inclusion
$\m{i}\subset \m{h}\times\m{v}$. Then the dimension of the face
$\m{F}_{i}$ is given by
\[\dim(\m{F}_\m{i}) = |\m{h}| + |\m{v}| - 1.\] Furthermore, assume that
$v \in \R_{\geq 0}\cdot\iw{S}{n}$ and that \[v = \sum_{(i,j) \in
\m{i}}c_{ij}(d_i e_i + d_j \cc{e}_j)\] for positive terms $c_{ij}$.
Then $\m{F}_{\m{i}}$ is the smallest face of $\R_{\geq
0}\cdot\iw{S}{n}$ containing $v$, i.e., $v\in \Int(\m{F}_{\m{i}})$.
\end{lemma}
{\it Proof: } Let $\m{h}$ and $\m{v}$ be as in the Lemma.  Assume
that $\m{h}$ is a proper subset of $\{1,\ldots,k\}$, and fix $j' \in
\m{v}$. Adding a horizontal line $i'\not\in\m{h}$ introduces a new
element $d_{i'} e_{i'} + d_{j'} \cc{e}_{j'}$ not in the vector space
spanned by $\langle d_i e_i + d_j \cc{e}_j:
(i,j)\in\m{h}\times\m{v}\rangle$. Let $\m{h}'=\m{h}\cup\{i'\}$, then
$\dim(\m{F}_{\m{h}'\times\m{v}}) \ge \dim(\m{F}_{\m{h}\times\m{v}})
+ 1$ by the above. A similar argument holds for the vertical lines.
Since $\dim(\m{F}_{\{1,\ldots,k\}\times\{1,\ldots,k\}}) = 2k-1$, it
follows that $\dim(\m{F}_\m{i}) = |\m{h}| + |\m{v}| - 1$. The second
part of the Lemma is a consequence of Lemma \ref{le-2} and Theorem
\ref{th-1}, since $\m{F}_{\m{i}}$ is the smallest face containing
each of the elements $d_i e_i + d_j \cc{e}_j$ for $(i,j)\in\m{i}$.
$\Box$\medskip

\begin{figure}[h]
  $$\includegraphics[scale = 1]{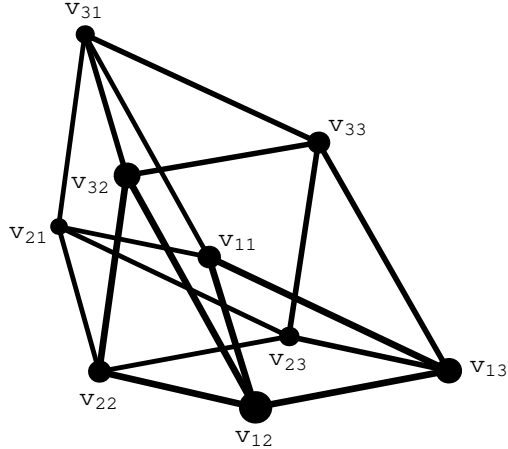}$$
  \caption{The cone $\R_{\geq
0}\cdot S_{n_1,n_2,n_3}$ in $\R^5$}\label{fig1}
\end{figure}

Using Lemma \ref{le-3} one counts the number of faces of $\R_{\geq
0}\cdot \iw{S}{n}$ using a simple combinatorial argument. The number
of faces $m_d$ of dimension $d$ with $0\leq d \leq 2k-1$ is given by
the formula\footnote{Count all choices of $p$ horizontal lines and
$q$ vertical lines with $p+q = d+1$.},\[m_{d} =
\binom{k}{1}\binom{k}{d} + \cdots + \binom{k}{d}\binom{k}{1}.\] In
particular the number of faces of $\R_{\geq 0}\cdot \iw{S}{n}$ of
codimension one is $m_{2k-2} = \binom{k}{k-1}\binom{k}{k} +
\binom{k}{k}\binom{k}{k-1} = 2k$ and faces of codimension $2$ are
$m_{2k-3} =
\binom{k}{k-2}\binom{k}{k}+\binom{k}{k-1}\binom{k}{k-1}+\binom{k}{k}\binom{k}{k-2}
= k(2k-1)$ and so on.
\begin{example}[Visualizing the structure of $\R_{\geq 0}\cdot \iw{S}{n}$]\label{ex-1}\end{example}
Using Theorem \ref{th-1} and Lemma \ref{le-3} one can visualize the
structure of the cones $\R_{\geq 0}\cdot \iw{S}{n}$, at least for
the minimal generators case mentioned in the Theorem. Figure
\ref{fig1}, shows one way to structure the cone for $k=3$. The
extreme rays, ${\mathrm v_{ij}} = \R_{\geq 0}\cdot (d_i e_i + d_j
\cc{e}_j)$ are shown as nodes and the two-dimensional faces of
$\R_{\geq 0}\cdot S_{n_1,n_2,n_3}$ are shown as edges connecting
their generating rays/vectors.

For the case, $k=3$, shown on Figure \ref{fig1} the number of faces
is counted as follows: $m_1 = 9$ (nodes on the Figure), $m_2 = 18$
(shown as edges), $m_3 = 15$ (appear as triangles and
quadrilaterals), $m_4 = 6$ (exclude a triangle from the Figure).

\section{Extending to affine toric varieties}
\label{se-3}

A maximal ideal $m$ in  $\C[\iw{S}{n}]$ determines, naturally, a
ring homomorphism $\C[\iw{S}{n}] \rightarrow \C~~[=\C[\iw{S}{n}]/m]$
and therefore induces a homomorphism between the monoids
$(\iw{S}{n},0,+)$ and $(\C,1,\cdot)$. Vice versa, a homomorphism
$u:(\iw{S}{n},0,+)\rightarrow (\C,1,\cdot)$ extends to a ring
homomorphism $\C[\iw{S}{n}] \rightarrow \C$, thereby determining a
maximal ideal $m = \ker(u)$ in $\C[\iw{S}{n}]$. The constructions
are bijective and it follows that the maximum ideals can be
identified with the set of all homomorphisms
$u:(\iw{S}{n},0,+)\rightarrow (\C,1,\cdot)$ denoted by
$\Hom(\iw{S}{n},\C)$ or $\Homsg(\iw{S}{n},\C)$ here. Define the
``conjugate'' $\cc{u}$ of an element $u\in \iw{M}{n}$ by
interchanging $e_i$ and $\cc{e}_i$ for $i=1,\ldots, k$, e.g.,
$\overline{e_i + \cc{e}_j} = e_j + \cc{e}_i$. This defines a
conjugate operator on the lattice such that $\cc{\iw{S}{n}} =
\iw{S}{n}$. Define an injection of the quotient space $\C^k/\T$ into
$\Hom(\iw{S}{n},\C)$ of the circle action in the minimal generators
case as follows: Take an element $(z_1,\ldots,z_k)$ from the orbit
of an element $z$ in $\C^k/\T$ and define, using the invariant
generators, $u(z)\in \Hom(\iw{S}{n},\C)$ by $u(z)(e_i + \cc{e}_i)=
z_i \cc{z}_i$ and $u(z)(d_i e_i + d_j \cc{e}_j) = z_i^{d_i}
\cc{z}_j^{d_j}$. The mapping $z\mapsto u(z)$ is
injective\footnote{If u(z) = u(w) then either both $z_i$ and $w_i$
are zero or $z^{d_i}_i w^{-d_i}_i$ is a unit constant $c =
\exp{\imath\theta}$ independent of $i$. Assuming that all the pairs
are nonzero, write $z_i w_i^{-1} = \exp(\imath\frac{\theta+2\pi
q_i}{d_i})$ for some integer $q_i$ and find integers
$s_1,\ldots,s_k$ so that $q_i + d_i s_i$ is independent of $i$,
i.e., it would be in $(q_1 + d_1 \Z)\cap (q_2 + d_2 \Z)\cap\cdots
\cap (q_k + d_k \Z)$ which is of the form $q + d_1 \cdots d_k \Z$
for some integer $q$. Finally, for $t= \exp(\imath\frac{\theta +
2\pi q}{d_1 \cdots d_k})$ it holds that $z_i = t^{n_i}w_i$ so $z$
and $w$ are in the same $\T$-orbit. The case when some of the pairs
$z_i$ and $w_i$ are zero also follows from this one.},
interestingly, it is injective independent of the assignments
$u(z)(e_i + \cc{e}_i)$. More generally, Hilbert basis resulting from
linear actions of compact groups separate orbits, see for example
\cite{Po}. But, as pointed out, only a simple argument is required
here to show this for the orbits of the circle action at hand.

\begin{theorem}\label{th-2}
Assume that the integers $n_1,\ldots,n_k$ generate a minimal Hilbert
basis. The injection of the orbit space $\C^k/\T$ into
$\Hom(\iw{S}{n},\C)$ consists of all $x$ in $\Hom(\iw{S}{n},\C)$
with $x(\cc{v}) = \overline{x(v)}$ for all $v\in \iw{S}{n}$ and
$x(e_i+\cc{e}_i) \geq 0$ for $i=1,\ldots,k$.
\end{theorem}
{\it Proof: } Denote by $\tau$ the element in $\Hom(\iw{S}{n},\C)$
defined by $0\mapsto 1$ and $v\mapsto 0$ if $v\neq 0$. For
$(z_1,\ldots,z_k)$ from the orbit of some element $z$ in $\C^k/\T$,
$x=u(z)$ satisfies the stated conditions on the generating set
$\{e_i+\cc{e}_i\} \cup \{d_i e_i + d_j \cc{e}_j\}$ of $\iw{S}{n}$
and therefore also on all of $\iw{S}{n}$. For the other direction:
Take $x\in \Hom(\iw{S}{n},\C)$ satisfying the conditions from the
Theorem. First if $x(e_i+\cc{e}_i) = 0$ for $i=1,\ldots,k$ then
$|x(d_i e_i+d_j\cc{e}_j)|^2 = x(
e_i+\cc{e}_i)^{d_i}x(e_j+\cc{e}_j)^{d_j} = 0$ for all $i,j$, so $x =
\tau$. Now assume that $x\neq \tau$, then for some $i$,
$x(e_i+\cc{e}_i)
> 0$. For simplicity, assume $x(e_1+\cc{e}_1) > 0$. Create a
representative $(w_1,\ldots,w_k)$ for $w\in\C^{k}/\T$ by $w_1 =
\sqrt{x(e_1+\cc{e}_1)}$ and for $i=2,\ldots,k$ let $w_i$ be any of
the $d_i$-th roots of $x(d_i e_i+d_1 \cc{e}_1)/w_1^{d_1}$. Then
$x(d_i e_i+d_1 \cc{e}_1) = w_i^{d_i} \cc{w_1}^{d_1}$.  It follows
that $x(d_i e_i+d_j \cc{e}_j) = x(d_i e_i+d_1 \cc{e}_1)x(d_1 e_1+d_j
\cc{e}_j)/x(d_1 e_1+d_1 \cc{e}_1) = x(d_i e_i+d_1
\cc{e}_1)\overline{x(d_j e_j+d_1 \cc{e}_1)}/x(e_1+ \cc{e}_1)^{d_1} =
w_i^{d_i} \cc{w}_j^{d_j}$. Finally, since $x(e_i+ \cc{e}_i) \geq 0$
and $x(e_i+ \cc{e}_i)^{d_i} = x(d_i e_i+ d_i \cc{e}_i) = (w_i
\cc{w}_i)^{d_i}$ it follows that $x(e_i+ \cc{e}_i) = w_i \cc{w}_i$.
In other words $x=u(w)$ is in the image of $u$ as required.
$\Box$\medskip

By identifying $\C^k/\T$ with its image under $u$ one may write
\[ \C^k/\T \subset
\Homsg(\iw{S}{n},\C).\] This identifies the orbit spaces as a nice
subset of the maximal ideals, or points, of the affine toric variety
$\Spec(\C[\iw{S}{n}])$, see Fulton \cite{Fu} section 1.3 for a
starting point. The description facilitates a study of the orbit
space using the geometry of toric varieties.

\begin{example}[The dual cone $\iw{S}{n}^\vee$]\label{ex-2}\end{example}
Here we describe the dual lattice of $\iw{M}{n}$, which is
identified with $\iw{M}{n}$ using the pairing defined below, and the
dual cone of $\R_{\geq 0}\cdot\iw{S}{n}$. Using the basis
$l_1,\ldots,l_k,\eta_1,\ldots,\eta_{k-1}$ from Lemma \ref{le-1},
calculate
\[v_{ij} = d_i e_i + d_j \cc{e}_j = \Big{\{}
\begin{array}{l}
d_i l_i + \cdots + d_j l_j - \eta_i - \cdots - \eta_{j-1}\mbox{ if
}i\leq j, \\ \eta_j + \cdots + \eta_{i-1} - d_{j+1} l_{j+1} - \cdots
- d_{i-1} l_{i-1} \mbox{ otherwise.}
\end{array}\]
Using the basis, these calculations and by identifying the lattice
$\iw{M}{n}$ with its dual lattice one obtains an inner product
$\Eval{}{}$ on $\iw{M}{n}$ satisfying
\[\begin{array}{ccc}
\Eval{v_{ij}}{l_m}=\Bigg{\{}
\begin{array}{cl}
+ d_m & \mbox{if }i\leq m \leq j,\\
- d_m & \mbox{if }j < m < i,\\
0 & \mbox{otherwise}
\end{array} &\mbox{and}&
\Eval{v_{ij}}{\eta_m}=\Bigg{\{}
\begin{array}{cl}
- 1 & \mbox{if }i\leq m < j,\\
+ 1 & \mbox{if }j \leq m < i,\\
0 & \mbox{otherwise.}
\end{array}
\end{array}
\]

An element $x = a_1 l_1 + \cdots a_k l_k + b_1 \eta_1 + \cdots +
b_{k-1} \eta_{k-1}$ is nonnegative on the semigroup $\iw{S}{n}$ if
$\Eval{v_{ij}}{x} \geq 0$ always, and, by the above,
\[\Eval{v_{ij}}{x} = a_i d_i + \ldots + a_j d_j - b_i - \cdots - b_{j-1}
\mbox{ if }i \leq j\]
\[\Eval{v_{ij}}{x} = - a_{j+1} d_{j+1} - \cdots - a_{i-1} d_{i-1} + b_j +
\cdots + b_{i-1} \mbox{ if } j < i.\] The set of all lattice points
that are nonnegative on $\iw{S}{n}$ under the pairing is denoted by
$\iw{S}{n}^\vee$, it spans a strongly convex rational polyhedral
cone $\R_{\geq 0} \cdot \iw{S}{n}^\vee$ dual to the cone $\R_{\geq
0} \cdot \iw{S}{n}$. So, by the above formula, if $x$ is in
$\R_{\geq 0} \cdot \iw{S}{n}^\vee$ then all $a_i\geq 0, b_j\geq 0$,
but this is not a sufficient condition for $x$ to belong to the dual
cone. The elements $x_1,\ldots,x_k$ and $y_1,\ldots,y_k$ given by
\[
x_1 = l_1 \mbox{ and } x_2 = l_2 + d_2\eta_1, \ldots, x_k = l_k +
d_{k}\eta_{k-1},\]
\[
y_1 = l_1 + d_1\eta_1, \ldots, y_{k-1} = l_{k-1} + d_{k-1}\eta_{k-1}
\mbox{ and } y_k = l_k\] are in $\iw{S}{n}^\vee$. It follows from
Lemma \ref{le-2} and the above that each of these $2k$ elements
spans a one-dimensional face of $\R_{\geq 0}\cdot\iw{S}{n}^\vee$.
The number of such rays is equal to the number of codimensional one
(facets) of the original cone $\R_{\geq 0}\cdot\iw{S}{n}$ which was
determined to be $m_{2k-2} = 2k$. Therefore these are all the rays.
\medskip

By denoting the standard basis for $\Z^{k^2}$ as $\{e_{ij}:~i,j =
1,\ldots,k\}$ define an epimorphism $F_k: \Z^{k^2}\rightarrow
\iw{M}{n}$ given by $e_{ij}\mapsto d_i e_i + d_j \cc{e}_j$ if $i\neq
j$ and $e_{ii}\mapsto e_i + \cc{e}_i$. This function is used
extensively in the remaining sections of the paper.
\medskip

\begin{example}[$F_3$]\label{ex-3}\end{example}
The map $F_3$ has matrix representation

$$F_3 = \begin{array}{cc}
          &
\begin{array}{ccccccccc}
  e_{21} & e_{32} & e_{11} & e_{22} & e_{33} & e_{31} & e_{12} & e_{23} & e_{13}
  \\ \\
\end{array} \\
\begin{array}{c}
  \eta_1 \\
  \eta_2 \\
  l_1 \\
  l_2 \\
  l_3 \\
\end{array} &
\left[
\begin{array}{ccccccccc}
  ~1~&~0~&~0~&~0~&~0~&~1~&-1~&~0~&-1~\\
  0 & 1 & 0 & 0 & 0 & 1 & 0 & -1 & -1 \\
  0 & 0 & 1 & 0 & 0 & 0 & d_1 & 0 & d_1 \\
  0 & 0 & 0 & 1 & 0 & -d_2 & d_2 & d_2 & d_2 \\
  0 & 0 & 0 & 0 & 1 & 0 & 0 & d_3 & d_3 \\
\end{array}\right]\end{array}$$ in the bases shown and its kernel is given by the image
of the matrix
$$
\begin{array}{cc}
\begin{array}{c}
  e_{21} \\ e_{32} \\ e_{11} \\ e_{22} \\ e_{33} \\ e_{31} \\ e_{12} \\ e_{23} \\ e_{13}
\end{array}&
\left[%
\begin{array}{cccc}
  1 & -1~ & 0 & -1~ \\
  1 & 0 & -1~ & -1~ \\
  0 & ~d_1 & 0 & d_1 \\
  -d_2~ & ~d_2 & ~d_2 & ~d_2 \\
  0 & 0 & ~d_3 & ~d_3 \\
  -1~ & 0 & 0 & 0 \\
  0 & -1~ & 0 & 0 \\
  0 & 0 & -1~ & 0 \\
  0 & 0 & 0 & -1~ \\
\end{array}%
\right] \\
\end{array}.$$
\medskip

Given a finite sequence of integers $n_1,\ldots, n_k$ generating a
minimal Hilbert basis, one may always extend it to a longer such
sequence by adding the integer $d = d_1 d_2 \cdots d_k$ repeatedly
to the sequence. The resulting cones $\R_{\geq 0}\cdot
S_{n_1,\ldots,n_k,d,\ldots,d}$ contain the original cone $\R_{\geq
0}\cdot \iw{S}{n}$ as one of its faces according to Theorem
\ref{th-1}. Similarly, the supporting lattice
$M_{n_1,\ldots,n_k,d,\ldots,d}$ contains $\iw{M}{n}$ by identifying
the $2k-1$ vectors $l_1, \ldots, l_k, \eta_1,\ldots,\eta_{k-1}$ from
Lemma \ref{le-1} for both lattices. The orthogonal projection $\pi:
M_{n_1,\ldots,n_k,d,\ldots,d} \rightarrow \iw{M}{n}; \sum t_i l_i +
\sum s_j \eta_j \mapsto \sum_{i=1}^{k} t_i l_i + \sum_{j=1}^{k-1}
s_j \eta_j$ is the identity map on $\iw{S}{n}$ and so maps
$S_{n_1,\ldots,n_k,d,\ldots,d}$ surjectively onto $\iw{S}{n}$. It
also induces an injection $\pi^{\ast}:\Homsg(\iw{S}{n},\C)
\hookrightarrow \Homsg(S_{n_1,\ldots,n_k,d,\ldots,d},\C)$.

\section{Poisson embeddings and the polyhedra}
\subsection{Faces resulting from squarefree generators}


The standard Poisson structure on $\R^{2k}$ extends to the invariant
polynomials $\C[\iw{S}{n}]$ as described in \cite{ASE1} and is given
there by the simple bracket
$$
\Pbr{X^a}{X^b} = -2\i\sum_1^k (a_i\cc{b}_i - \cc{a}_i b_i)X^{a + b -
l_i}
$$ on monomials.
This formula does not appear to represent the symmetry seen on Figure
\ref{fig1} from Example \ref{ex-1} very well, since, on the one
hand, it favors the extremal rays generated by $l_i$ $(v_{ii})$ over
$v_{ij},~(i\neq j)$. On the other, each of the extremal rays
$v_{ij}$ is attached to a similar structure of faces: referring to
Figure \ref{fig1} for $k=3$, this structure is always four lines
(2-dim faces), two triangles and four quadrilaterals (3-dim faces),
and four 4-dim facets each obtained by removing one of the other
four triangles.

A family of Poisson brackets, that includes the standard Poisson
algebra and is based on the structure of the polyhedral cone
$\R_{\geq 0} \iw{S}{n}$ is defined in the following lemma. It is
assumed that the integers $n_1,\ldots, n_k$ generate a minimal
Hilbert basis and as before: $d_i =
\gcd(n_1,\ldots,n_{i-1},n_{i+1},\ldots,n_k)$. The lemma connects
squarefree monomials (e.g., see \cite{St} for a discussion of
Stanley-Reisner ideals) from the Hilbert basis of the $\T$
invariants and Poisson structures on the $\T$ orbit space in a
simple way.

\begin{lemma}\label{le-4}
Let $\epsilon = (\epsilon_{ij})$ be a real $k\times k$ matrix
satisfying $\epsilon_{ij} = 0$ if $d_i \neq d_j$. The bilinear
antisymmetric bracket $\Pbr{~}{}_\epsilon$ on $\C[\iw{S}{n}]$
determined on monomials by
$$ \Pbr{X^a}{X^b}_\epsilon = -2\i\sum_{ij} \epsilon_{ij}(a_i\cc{b}_j
- \cc{a}_j b_i)X^{a + b - e_i - \cc{e}_j}
$$ is a Poisson bracket on $\C[\iw{S}{n}]$.
\end{lemma}
{\it Proof: } Note first that the formula $d_i = d_j$ always implies
that $X^{e_i + \cc{e}_j}$ is an invariant since either $i=j$ or
$d_i=d_j=1$ if $d_i$ and $d_j$ are equal. If for $a,b \in \iw{S}{n}$
and $i,j$ with $d_i = d_j$ it holds that $a+b-e_i-\cc{e}_j
\not\in\iw{S}{n}$, then $a_i + b_i - 1 < 0$ or $\cc{a}_j + \cc{b}_j
- 1 < 0$, but since $a_i, b_i, \cc{a}_j$ and $\cc{b}_j$ are all
nonnegative it follows that $a_i = b_i = 0$ or $\cc{a}_j = \cc{b}_j
= 0$ and therefore $(a_i\cc{b}_j - \cc{a}_j b_i)X^{a + b - e_i -
\cc{e}_j} = 0$. This shows that the bracket $\Pbr{~}{~}_\epsilon$
maps $\C[\iw{S}{n}]\times\C[\iw{S}{n}]$ into $\C[\iw{S}{n}]$.
Leibniz identity\footnote{$\Pbr{f g}{h} = f\Pbr{g}{h} +
\Pbr{f}{h}g$} follows from verifying the formula $\Pbr{X^a}{X^b
X^c}_\epsilon = \Pbr{X^a}{X^b}_\epsilon X^c +
\Pbr{X^a}{X^c}_\epsilon X^b$ for $a,b,c\in \iw{S}{n}$ directly.
Formally, write $\Pbr{X^a}{X^b}_\epsilon =$ $$ -2\i\sum_{ij}
\epsilon_{ij} (\p{X^a}{z_i}\p{X^b}{\cc{z}_j} -
\p{X^a}{\cc{z}_j}\p{X^b}{z_i}) =
-2\i\sum_{ij}\epsilon_{ij}\DD{z_i}{\cc{z}_j}(dX^a \wedge dX^b).$$
Using this formula one verifies Jacobi
identity\footnote{$\Pbr{\Pbr{f}{g}}{h} + \Pbr{\Pbr{g}{h}}{f} +
\Pbr{\Pbr{h}{f}}{g} = 0$} by extending the bracket to all
polynomials in the variables $z_i$ and $\cc{z}_j$, i.e.,
$\Pbr{z_i}{\cc{z}_j} = -2\i \epsilon_{ij}$. The identity now follows
from antisymmetry and Leibniz identity. Consequently
$(\C[\iw{S}{n}],\Pbr{}{}_\epsilon)$ is a Poisson algebra. $\Box$

Referring the reader to section 6 in \cite{ASE1} the Poisson
bivector above is converted into real coordinates $x,y$ on
$\R^{2k}$, satisfying $x_i + \i y_i = z_i$ and $x_i - \i y_i =
\cc{z}_i$, by replacing $\D{z_i}$ with $\frac{1}{2}(\D{x_i} - \i
\D{y_i})$ and $\D{\cc{z}_i}$ with $\frac{1}{2}(\D{x_i} + \i
\D{y_i})$. Assuming that $\epsilon$ is symmetric this results in the
coordinate change
$$-2\i\sum_{ij}\epsilon_{ij}\DD{z_i}{\cc{z}_j}
= \sum_{ij}\epsilon_{ij}\DD{x_i}{y_j}.$$ In particular, the right
hand side bivector is a true real valued, as appose to complex
valued, Poisson bivector on $R^{2k}$ - as a result of the symmetry
of $\epsilon$. Its rank on $\R^{2k}$ is given by the the rank of the
matrix $\left(
\begin{array}{cc}
0 & \epsilon \\
-\epsilon & 0 \\
\end{array}
\right)$.

The epimorphism $F_k: \Z^{k^2}\rightarrow \iw{M}{n}$, defined in
section \ref{se-3}, induces an algebra morphism
$\C[X_{ij}]\rightarrow \C[\iw{S}{n}]$ by $F_k (X_{ij}) =
X^{F_k(e_{ij})}$. A Poisson bivector on $\C[X_{ij}]$ is considered
real valued, see \cite{ASE1}, if it contains no imaginary part after
being transformed into real coordinates according to the mappings:
\begin{itemize}
\item $X_{ij} \mapsto R_{ij} + \i I_{ij}$ and $X_{ji}
\mapsto R_{ij} - \i I_{ij}$ for $i<j$,
\item $\D{Z_{ij}} \mapsto
\frac{1}{2}(\D{R_{ij}} -\i\D{I_{ij}})$ and $\D{Z_{ji}} \mapsto
\frac{1}{2}(\D{R_{ij}} + \i\D{I_{ij}})$ for $i<j$.
\end{itemize}
Such a real bivector determines a Poisson algebra in the $k^2$
variables $X_{ii}, R_{i<j}$ and $I_{i<j}$ over $\R$ and determines a
Poisson structure on $\R^{k^2}$.

\begin{theorem}\label{th-3}
A Poisson structure of the form $\Pbr{~}{~}_{\epsilon}$ on
$\C[\iw{S}{n}]$, with $\epsilon$ symmetric and $\epsilon_{ij} = 0$
unless $d_i = d_j = 1$, is $F_k$ related\footnote{$F_k(\Pbr{f}{g}) =
\Pbr{F_k(f)}{F_k(g)}_\epsilon$} to a real Poisson structure on
$\C[X_{ij}]$.
\end{theorem}
{\it Proof:} Define a bilinear bracket $\Pbr{~}{~}$ on $\C[X_{ij}]$
using the formula $$\Pbr{X_{pq}}{X_{st}} = -2\i(\epsilon_{pt} X_{sq}
- \epsilon_{sq} X_{pt})$$ and by extending it to all the polynomials
by way of Leibniz identity and bilinearity. It is $F_k$ related to
$\Pbr{~}{~}_{\epsilon}$ since $F_k(\Pbr{X_{pq}}{X_{st}}) =
-2\i(\epsilon_{pt} X^{F_k(e_{sq})} - \epsilon_{sq} X^{F_k(e_{pt})})
= \Pbr{F_k(e_{pq})}{F_k(e_{st})}_{\epsilon}$, as a result of the
condition $\epsilon_{ij} = 0$ if not both $d_i$ and $d_j$ are equal
to one.  The new bracket satisfies Jacobi identity because its
Jacobiator $\m{J}$\footnote{$\m{J}(A,B,C) = \Pbr{\Pbr{A}{B}}{C} +
\Pbr{\Pbr{B}{C}}{A} + \Pbr{\Pbr{C}{A}}{B}$} maps triplets $A,B,C$
from the set of indeterminants $\{X_{ij}\}$ into linear polynomials
which again map to zero under $F_k$ since $F_k \circ \m{J} = 0$ by
Jacobi identity for the original bracket $\Pbr{~}{~}_\epsilon$. The
only linear polynomial that maps to zero under $F_k$ is zero itself
so it follows that Jacobi identity is also satisfied for the derived
bracket $\Pbr{~}{~}$. The bracket $\Pbr{~}{~}$ is real since it is
the unique lift under $F_k$ of the real Poisson bracket
$\Pbr{~}{~}_\epsilon$ to a linear Poisson structure on $\C[X_{ij}]$.
$\Box$

\subsection{Intertwining other faces}

Theorem \ref{th-3} connects the $\T$ invariant Poisson structure
$\Pbr{}{}_\epsilon$ on $\R^{2k}$ and the special face
$\m{F}_{\m{l}\times\m{l}}$ where $\m{l} = \{~i~:~d_i = 1~\}$ of the
polyhedral cone $\R_{\geq 0}\cdot\iw{S}{n}$. In order to accommodate
the other faces also, below consider a generalization denoted by
$\Pbr{}{}_\epsilon^\delta$ and given by:
$$ \Pbr{X^a}{X^b}_\epsilon^\delta = -2\i\sum_{ij} \epsilon_{ij}(a_i\cc{b}_j
- \cc{a}_j b_i)X^{a + b + \delta_{ij}F_k(e_{ij})}$$ where $\delta$
is taken to be an integer matrix. The bracket is extended to all of
$\C[z,\cc{z}]$ by the formulas $\Pbr{z_i}{\cc{z}_j}_\epsilon^\delta
= -2\i\epsilon_{ij} (z_i \cc{z}_j) X^{\delta_{ij}F_k(e_{ij})}$ and
zero on other pairs and via bilinearity, antisymmetry and Leibniz
identity. Consequently, it is required that $\delta_{ij} \geq -1$ if
$d_i = d_j$ and $\delta_{ij} \geq 0$ otherwise. Requiring both
$\epsilon$ and $\delta$ to be symmetric and then converting the
resulting bivector to real coordinates results in a real bracket as
follows $\sum_{ij} \Pbr{z_i}{\cc{z}_j}_\epsilon^\delta
\DD{z_i}{\cc{z}_j} =$
$$ \frac{1}{2}\sum_{i<j}\Re\Pbr{z_i}{\cc{z}_j}_\epsilon^\delta ~
(\DD{x_i}{x_j}+\DD{y_i}{y_j})-\frac{1}{2}\sum_{ij}\Im\Pbr{z_i}{\cc{z}_j}_\epsilon^\delta
~ \DD{x_i}{y_j}.$$ To establish what conditions are needed for
Jacobi identity to hold, first note that the Jacobiator for the
bracket is always zero when applied to triplets of the form
$(z_p,z_q,z_r)$ and $(\cc{z}_p,\cc{z}_q,\cc{z}_r)$. For the mixed
triplets $(z_p,z_q,\cc{z}_r)$ and $(\cc{z}_p,\cc{z}_q,z_r)$ it holds
that $\m{J}(z_p,z_q,\cc{z}_r) =
\overline{\m{J}(\cc{z}_p,\cc{z}_q,z_r)} =$ $$ -4
\epsilon_{pr}\epsilon_{qr} (\delta_{pr} d_r^p - \delta_{qr}
d_r^q)z_p z_q \cc{z}_r X^{\delta_{pr}
F_k(e_{pr})+\delta_{qr}F_k(e_{qr})}$$ where the notation $d_i^j$ is
used to denote $d_i$ if $i\neq j$ and $d_i^i = 1$. The results are
summarized as follows:
\begin{corollary}\label{cr-1}
The bracket $\Pbr{~}{}_\epsilon^\delta$ is Poisson on the polynomial
algebra $\R[x,y]$ in $2k$ variables if: $\epsilon$ and $\delta$ are
symmetric matrices, $\epsilon$ has real coefficients, $\delta$ has
integer coefficients satisfying $\delta_{ij}\geq -1$ when $d_i =
d_j$ and otherwise $\delta_{ij} \geq 0$, and for each triplet
$p,q,r$ the equation $\epsilon_{pr}\epsilon_{qr} (\delta_{pr} d_r^p
- \delta_{qr} d_r^q) = 0$ holds. Furthermore, in this case, if $f$
and $g$ are $\T$ invariant polynomials then so is
$\Pbr{f}{g}_\epsilon^\delta$.
\end{corollary}

\subsection{Explicit structures and lifts}

\begin{example}[$\delta_{ii} = d_i$ and $\delta_{ij}=1$ if $i\neq j$]
\label{ex-5}\end{example} A simple way to have
$\Pbr{~}{}_\epsilon^\delta$ satisfy Jacobi identity is to fix
$\delta$ by defining $\delta_{ii} = d_i$ and $\delta_{ij}=1$ if
$i\neq j$, in this case the equations $\epsilon_{pr}\epsilon_{qr}
(\delta_{pr} d_r^p - \delta_{qr} d_r^q) = 0$ are always satisfied.
The resulting Poisson bracket may be written out as
$\Pbr{z_i}{\cc{z}_j}_\epsilon^\delta = -2\i\epsilon_{ij}~
z^{d_i+1}_i \cc{z}^{d_j+1}_j$.

\begin{example}[$\delta = 0$]\label{ex-6}\end{example}
The bracket $\Pbr{z_i}{\cc{z}_j}_\epsilon^\delta$ satisfies Jacobi
identity if $\delta$ is the zero matrix. Letting $\delta=0$ results
in a bracket given by $\Pbr{z_i}{\cc{z}_j}_\epsilon^0 = -2\i
\epsilon_{ij}~z_i \cc{z}_j$. Now lift $\Pbr{~}{~}_\epsilon^0$ to
$\C[X_{ij}]$ as follows: Define a bilinear bracket $\Pbr{~}{~}$ on
$\C[X_{ij}]$ using the formulas
$$\Pbr{X_{pq}}{X_{st}} = -2\i(\epsilon_{pt} d_p^q d_t^s  -
\epsilon_{sq}d_q^p d_s^t)X_{pq} X_{st}$$ and extend it to all the
polynomials using Leibniz identity. The new bracket is defined in
such a way that it is $F_k$ related to $\Pbr{~}{~}_{\epsilon}^0$.
The new bracket also satisfies Jacobi identity: Calculating
$\Pbr{\Pbr{X_{pq}}{X_{st}}}{X_{ij}}$ results in
$\Pbr{\Pbr{X_{pq}}{X_{st}}}{X_{ij}} = -4 E_{pt}^{sq} (E_{pj}^{iq} +
E_{sj}^{it})X_{pq}X_{st}X_{ij}$ where $E_{ad}^{cb} = \epsilon_{ad}
d_a^b d_d^c  - \epsilon_{cb}d_b^a d_c^d$ and Jacobi identity for
$\Pbr{~}{~}$ now follows from calculating the other parts of the
Jacobiator $\m{J}(X_{pq},X_{st},X_{ij})$ and using $E_{ad}^{cb} = -
E_{cb}^{ad}$ to cancel terms - or, even simpler, by using that $F_k
\circ \m{J}=0$ . When real coordinates are introduced on
$\C[X_{ij}]$, see discussion before Theorem \ref{th-3}, the
resulting conjugate operator satisfies $\cc{X}_{ij} = X_{ji}$ and
the bracket therefore satisfies $\overline{\Pbr{X_{pq}}{X_{st}}} =
\Pbr{\cc{X}_{pq}}{\cc{X}_{st}}$ assuming that $\epsilon$ is
symmetric. This condition guaranties that the bracket is real valued
and as such restricts to a $k^2$ dimensional Poisson algebra on
$\R[X_{ii},R_{i<j},I_{i<j}]$.

\begin{example}[$\delta_{ij} = -1$ if $d_i = d_j = 1$ and $\delta_{ij}=0$ otherwise]
\label{ex-7}\end{example} If $\delta$ is fixed as: $\delta_{ij} =
-1$ if $d_i = d_j = 1$ and $\delta_{ij} = 0$ otherwise, then the
formulas $\epsilon_{pr}\epsilon_{qr} (\delta_{pr} d_r^p -
\delta_{qr} d_r^q) = 0$ are satisfied by requiring additionally that
$\epsilon_{ij} = 0$ whenever exactly one of $d_i$ and $d_j$ is equal
to one. The resulting Poisson bracket is determined by $$
  \Pbr{z_i}{\cc{z}_j}_\epsilon^\delta = -2\i \epsilon_{ij} \mbox{ if } d_i =
d_j = 1,$$
  $$\Pbr{z_i}{\cc{z}_j}_\epsilon^\delta = -2\i
\epsilon_{ij}~z_i \cc{z}_j \mbox{ if } d_i\neq 1\mbox{ and }d_j \neq
1.$$ It may be lifted to an $F_k$ related real bracket $\Pbr{~}{~}$
on $\C[X_{ij}]$ determined by $\Pbr{X_{pq}}{X_{st}} =
-2\i(\epsilon_{pt} P_{pt}^{sq} -\epsilon_{sq} P_{sq}^{pt})$ where
$P_{pt}^{sq} = X_{sq}$ if $d_p=d_t=1$ and $P_{pt}^{sq} = d_p^q d_t^s
X_{pq}X_{st}$ otherwise (similarly $P_{sq}^{pt} = X_{pt}$ if
$d_s=d_q=1$ and $P_{sq}^{pt} = d_s^t d_q^p X_{pq}X_{st}$ otherwise).
This formula may also be written
$$\Pbr{X_{pq}}{X_{st}} =
-2\i(\epsilon_{pt} d_p^q d_t^s {\Big\{}\begin{array}{c}
X_{sq}~^\mathrm{I}\\
X_{pq}X_{st}~^\mathrm{II}\\
\end{array}
-\epsilon_{sq}d_q^p d_s^t{\Big\{}
\begin{array}{c}
X_{pt}~^\mathrm{i}\\
X_{pq}X_{st}~^\mathrm{ii} \\
\end{array})$$ where the monomials shown are selected based on the following schema:
polynomial I is used when $d_p=d_t=1$ and II is used otherwise, also
polynomial i is used when $d_s=d_q=1$ otherwise polynomial ii is
used to complete the formula. As in the previous example, this
bracket is seen to be real valued since when real coordinates are
introduces the resulting conjugate operator satisfies
$\overline{\Pbr{X_{pq}}{X_{st}}}=\Pbr{\cc{X}_{pq}}{\cc{X}_{st}}$. In
order to prove that the new bracket satisfies Jacobi identity, write
$\Pbr{X_{pq}}{X_{st}} = \Pbr{X_{pq}}{X_{st}}_A +
\Pbr{X_{pq}}{X_{st}}_B$ where $\Pbr{~}{~}_A$ and $\Pbr{~}{~}_B$ are
given by $\Pbr{X_{pq}}{X_{st}}_A = -2\i(\epsilon^A_{pt}X_{sq} -
\epsilon^A_{sq}X_{pt})$ and $\Pbr{X_{pq}}{X_{st}}_B =
-2\i(\epsilon^B_{pt} d_p^q d_t^s - \epsilon^B_{sq}d_q^p d_s^t
)X_{pq}X_{st}$. Here $\epsilon^A_{ab}$ is zero unless $d_a = d_b =
1$ in which case $\epsilon^A_{ab} = \epsilon_{ab}$, similarly
$\epsilon^B_{ab}$ is zero unless $d_a \neq 1$ and $d_b \neq 1$ and
then $\epsilon^B_{ab} = \epsilon_{ab}$. As a result
$\epsilon^A_{ij}\epsilon^B_{it} = \epsilon^A_{ij}\epsilon^B_{sj}
=0$. The bracket $\Pbr{~}{~}_A$ is the one used in the proof of
Theorem \ref{th-3} so it already satisfies Jacobi identity. The
bracket $\Pbr{~}{~}_B$ is discussed in Example \ref{ex-6} and is
shown there to satisfy Jacobi identity. It follows that the
Jacobiator $\m{J}$ for $\Pbr{~}{~}$ satisfies
$\m{J}(X_{pq},X_{st},X_{ij}) =$
$$\Pbr{\Pbr{X_{pq}}{X_{st}}_B}{X_{ij}}_A+
\Pbr{\Pbr{X_{st}}{X_{ij}}_B}{X_{pq}}_A+
\Pbr{\Pbr{X_{ij}}{X_{pq}}_B}{X_{st}}_A+$$
$$\Pbr{\Pbr{X_{pq}}{X_{st}}_A}{X_{ij}}_B+
\Pbr{\Pbr{X_{st}}{X_{ij}}_A}{X_{pq}}_B+
\Pbr{\Pbr{X_{ij}}{X_{pq}}_A}{X_{st}}_B.$$ Using
$\Pbr{\Pbr{X_{pq}}{X_{st}}_B}{X_{ij}}_A =$ $$-4(
\epsilon^A_{sj}\epsilon^B_{pt} d_p^q d_t^s X_{pq} X_{it}
+\epsilon^A_{it}\epsilon^B_{sq} d_s^t d_q^p X_{sj} X_{pq}$$
$$-\epsilon^A_{pj}\epsilon^B_{sq} d_s^t d_q^p X_{iq} X_{st}
-\epsilon^A_{iq}\epsilon^B_{pt} d_p^q d_t^s X_{st} X_{pj})$$ and
$\Pbr{\Pbr{X_{pq}}{X_{st}}_A}{X_{ij}}_B =$ $$-4( \epsilon^A_{pt}(
\epsilon^B_{sj} d_s^t d_j^i - \epsilon^B_{iq} d_i^j
d_q^p)X_{sq}X_{ij}$$
$$-\epsilon^A_{sq}( \epsilon^B_{pj} d_p^q d_j^i
- \epsilon^B_{it} d_i^j d_t^s)X_{pt}X_{ij})$$ and similar formulas
for the other terms in the expression for
$\m{J}(X_{pq},X_{st},X_{ij})$ one concludes that all the terms
cancel when added together. It follows that the Jacobi identity for
$\Pbr{~}{~}$ is satisfied.

\subsection{Intertwined Poisson structures}
The real Poisson structure $\Pbr{~}{~}$ on $\C[X_{ij}]$ defined in
Example \ref{ex-7} will be denoted below as
$\Pbr{~}{~}^\m{F}_\epsilon$ where $\m{F} = \m{F}_{\m{l}\times\m{l}}$
denotes the face of $\R_{\geq 0}\cdot\iw{S}{n}$ indexed by $\m{l} =
\{~i~:~d_i = 1\}$. The same notation $\Pbr{~}{~}^\m{F}_\epsilon$
will be used to denote the real Poisson structure
$\Pbr{~}{~}_{\epsilon}^{\delta}$ on $\C[z,\cc{z}]$ for $\delta_{ij}
= -1$ if $d_i = d_j = 1$ and $\delta_{ij}=0$ otherwise. According to
the Example using this notation the brackets are $F_k$ related by
$$F_k(\Pbr{f}{g}^\m{F}_\epsilon) =
\Pbr{F_k(f)}{F_k(g)}^\m{F}_\epsilon$$ for polynomials $f,g$ in
$\C[X_{ij}]$. It is required that $\epsilon$ be symmetric and
$\epsilon_{ij}$ is zero if exactly one of $d_i$ and $d_j$ is one. In
order to simplify the notation further for $\epsilon = \Id$, the
identity matrix, the notation $\Pbr{~}{~}^\m{F}_\Id$ is reduced to
$\Pbr{~}{~}^\m{F}$ on both spaces. The bracket $\Pbr{~}{~}^\m{F}$ is
specified by the bivector $\Pi^{\m{F}}$, satisfying
$\Pbr{f}{g}^\m{F} = \Pi^{\m{F}}(df\wedge dg)$ under the usual
pairing of bivectors and 2-forms and given by $$\Pi^\m{F} =
-2\i\sum_{d_i=1}\DD{z_i}{\cc{z}_i}-2\i\sum_{d_i\neq
1}z_i\cc{z}_i\DD{z_i}{\cc{z}_i}.$$ The more general notation
$\Pi^\m{F}_{\epsilon}$ is used to denote the bivector dual to the
Poisson structure $\Pbr{~}{~}^\m{F}_\epsilon$ on $\C[X_{ij}]$ as
well as the Poisson bivector dual to the real bracket
$\Pbr{~}{~}^\delta_\epsilon$ on $\C[z,\cc{z}]$ where $\epsilon$ is
symmetric, $\delta_{ij} = -1$ if $d_i = d_j = 1$ and $\delta_{ij}=0$
otherwise, and $\epsilon_{ij} = 0$ if exactly one of $d_i$ and $d_j$
is equal to one. The following theorem has been proven.

\begin{theorem}\label{th-4}
Assume that the positive integers $n_1,\ldots,n_k$ generate a
minimal Hilbert basis with respect to the Hamiltonian circle action
$\T\times \R^{2k} \rightarrow \R^{2k}$ with weights
$n_1,\ldots,n_k$. Let each of $d_i$ be the greatest common divisor
of all the weights except $n_i$ and let $\m{F}$ be the face of the
polyhedral cone $\R_{\geq 0}\cdot\iw{S}{n}$ given by $\m{F} =
\m{F}_{\m{l}\times\m{l}}$ where $\m{l} = \{~i~:~d_i = 1 \}$. Then
the Hilbert map $F_k^{\ast} = (F_k(X_{ij}))_{ij}$ restricted to
$\R^{2k} \rightarrow \R^{k^2}$ projects to a Poisson embedding of
the orbit space $$(\R^{2k}/\T,\Pi^{\m{F}}_\epsilon) \hookrightarrow
(\R^{k^2},\Pi^{\m{F}}_\epsilon),$$ for any symmetric real $k\times
k$ matrix $\epsilon$ satisfying $\epsilon_{ij} = 0$ if $i\in\m{l}$
and $j\not\in\m{l}$.
\end{theorem}

About the Hilbert embedding: The term, ``embedding", is used since
the orbit space $\R^{2k}/\T$ may be assigned a smooth structure
$\cInf(\R^{2k}/\T)$ of smooth $\T$ invariant functions on $\R^{2k}$.
From a theorem by Schwarz in \cite{Sz} it follows that
$F_k^{\ast\ast}\cInf(\R^{k^2}) = \cInf(\R^{2k}/\T)$ and according to
Mather in \cite{Mr} the Hilbert map is a proper embedding. For a
more complete discussion of how to transfer differential geometry
methods to the orbit spaces see \cite{Pf}. Note also that in the
discussion before Theorem \ref{th-2} it is shown directly that the
Hilbert map is injective on the orbit space since it is essentially
the mapping $z\rightarrow u(z)$ from Theorem \ref{th-2} even if the
target space is different in that case.

The rank of the Poisson structure $\Pi^{\m{F}}$ on $\R^{2k}$ is
always between $2|\m{l}|$ and $2k$. Explicitly, the rank of the
Poisson structure at a point $z\in \C^k$ is just twice the number of
indexes $i$ satisfying $z_i \neq 0$ or $i\in\m{l}$.

The discussion above and the theorem can be extended to faces $\m{T}
= \m{F}_{\m{h}\times\m{h}}$ where $\m{h}$ is any subset of $\m{l}$.
The resulting bracket on $\R^{2k}$ is determined, similarly as
before, by
$$
  \Pbr{z_i}{\cc{z}_j}_\epsilon^\m{T} = -2\i \epsilon_{ij} \mbox{ if } i,j\in\m{h},$$
  $$\Pbr{z_i}{\cc{z}_j}_\epsilon^\m{T} = -2\i
\epsilon_{ij}~z_i \cc{z}_j \mbox{ otherwise }$$ and
$\epsilon_{ij}=0$ if exactly one of $i$ or $j$ is in $\m{h}$. This
Poisson structure then lifts under the Hilbert embedding to the
Poisson structure on $\C[X_{ij}]$ determined by:
$$\Pbr{X_{pq}}{X_{st}}_\epsilon^\m{T} =
-2\i(\epsilon_{pt} d_p^q d_t^s {\Big\{}\begin{array}{c}
X_{sq}~^\mathrm{I}\\
X_{pq}X_{st}~^\mathrm{II}\\
\end{array}
-\epsilon_{sq}d_q^p d_s^t{\Big\{}
\begin{array}{c}
X_{pt}~^\mathrm{i}\\
X_{pq}X_{st}~^\mathrm{ii} \\
\end{array})$$ where, similarly as before, polynomial I is used when $p,t\in\m{h}$ and II is used
otherwise, and polynomial i is used when $s,q\in{h}$ otherwise
polynomial ii is used.

In particular, if $\m{h}=\emptyset$, so $\m{F}_{\m{h}\times\m{h}} =
0$, one obtains the bracket $\Pbr{~}{~}_{\epsilon}^0$ from Example
\ref{ex-6} which is independent of the weights $n_1,\ldots,n_k$.
Even when the weights do not generate a minimal Hilbert basis, as
required in the above, the structure $\Pbr{~}{~}_{\epsilon}^0$ lifts
to the target space of the Hilbert embedding to a (product) Poisson
bracket related to $\Pi^0$ under the embedding.

\subsection{Additional squarefree generators}

Interestingly, the lifted Poisson structure $\Pi_{\epsilon}^{\m{T}}$
on $\C[X_{ij}]$ for $\m{T} = \m{F}_{\m{h}\times\m{h}}$ where $\m{h}$
is a subset of the indexes $i$ satisfying $d_i=1$ intertwines the
simple product Poisson structure from Example \ref{ex-6} and the
Poisson structure from Theorem \ref{th-3}. To be exact, if $A$ and
$B$ references the two Poisson structures this means that
$$~\footnote{$\circlearrowleft\Pbr{\Pbr{a}{b}_A}{c}_B= \Pbr{\Pbr{a}{b}_A}{c}_B + \Pbr{\Pbr{b}{c}_A}{a}_B
+ \Pbr{\Pbr{c}{a}_A}{b}_B$}\circlearrowleft\Pbr{\Pbr{a}{b}_A}{c}_B +
\circlearrowleft\Pbr{\Pbr{a}{b}_B}{c}_A=0$$ as seen in the
calculations in Example \ref{ex-7}.

Let $n_1,\ldots,n_k$ be positive integer weights generating a
minimal Hilbert basis and let $d_1,\ldots,d_k$ be as before. Define
$d$ to be the product $d=d_1\cdots d_k$. Let $\T$ act on $\C^k$ with
weights $n_1,\ldots,n_k$ and $\T$ act on $\C^{2k}$ with weights
$n_1,\ldots,n_k,d,\dots,d$. The action on $\C^{2k}$ also generates a
minimal Hilbert basis and to accommodate both actions the sequence
$d_1,\ldots,d_k$ is extended by letting $d_{k+1} = \cdots = d_{2k} =
1$. A Hilbert basis for the invariants of the action on $\C^{2k}$
contains $z_1 \cc{z}_1, \ldots, z_{2k}\cc{z}_{2k}$ and $z^{d_i}_i
\cc{z}^{d_j}_j$ with $i\neq j$ both from $\{1,\ldots,2k\}$; the
action has at least $k+k^2$ squarefree generators even if the
original action on $\C^{k}$ may only produce $k$ squarefree
generators. There are many ways to embed the space $\C^k$ into
$\C^{2k}$ in ways compatible with the $\T$ actions. Consider the
following embedding
$$\C^{k}\hookrightarrow \C^{2k};(z_1,\ldots,z_k) \mapsto
(z_1,\ldots,z_k,z^{d_1}_1,\ldots,z^{d_k}_k).$$ The embedding is $\T$
compatible so it projects to the orbit spaces. The real Poisson
structure $\Pi^{\m{T}}$ on $\C^{2k}$ derived from the cone $\R_{\geq
0} S_{n_1,\ldots,n_k,d,\ldots,d}$ by taking $\m{T}$ to be the face
$\m{T} = \m{F}_{\m{w}\times\m{w}}$ with $\m{w} = \{k+1,\ldots,2k\}$
has rank everywhere at least $2k$ and is determined by:
$$\Pbr{z_i}{\cc{z}_i}^{\m{T}} = -2\i z_i \cc{z}_i\mbox{ and
}\Pbr{z_{k+i}}{\cc{z}_{k+i}}^{\m{T}} = -2\i\mbox{ for
}i=1,\ldots,k.$$ Using the Hilbert embedding $F_{2k}$ one obtains:
$$\frac{\C^k}{\T} \hookrightarrow (\frac{\C^{2k}}{\T},\Pi^{\m{w}})
\begin{array}{c}F_{2k}\\ ~\hookrightarrow~\\~\end{array}
(\R^{4k^2},\Pi^{\m{w}}).$$

It is interesting to compare the above with the Poisson embedding
dimension, defined by Davis in \cite{BEN1} as the smallest possible
dimension of the target space on which there exists a Poisson
structure that pullback to the usual\footnote{That is determined by
$\Pbr{x_i}{y_i} = 1$.} Poisson structure on $\C^k$ under the Hilbert
embedding. In the previous sections, the nondegeneracy, i.e., full
rank everywhere ($2k$), condition has been relaxed somewhat in order
to allow Poisson embeddings of minimum embedding dimension. In this
last section the embedding problem is yet again modified, requiring
the target spaces to be of higher, but finite, dimension.

\end{document}